\documentclass[lettersize,journal]{IEEEtran}
\usepackage{amsmath,amsfonts}
\usepackage{algorithmic}
\usepackage{algorithm}
\usepackage{array}
\usepackage[caption=false,font=normalsize,labelfont=sf,textfont=sf]{subfig}
\usepackage{textcomp}
\usepackage{stfloats}
\usepackage{url}
\usepackage{verbatim}
\usepackage{graphicx}
\usepackage{cite}
\usepackage{multirow}
\usepackage{bm}
\usepackage{color}
\usepackage{epstopdf}
\graphicspath{{Fig/}}
\allowdisplaybreaks
\def\r{\mathbf{r}}
\def\y{\mathbf{y}}
\def\z{\mathbf{z}}
\def\D{\mathbf{D}}
\def\S{\mathbf{S}}

\def\mE{{\mathcal E}}

\def\p{\partial}
\def\bfp{\mathbf{p}}
\def\f{\frac}
\def\Om{\Omega}
\def\na{\nabla}
\def\bmsigma{\bm{\sigma}}
\def\argmin{\mathop{\arg\min}}

\newcommand{\sgn}[1]{\mathrm{sgn}(#1)}

\def\BibTeX{{\rm B\kern-.05em{\sc i\kern-.025em b}\kern-.08em
    T\kern-.1667em\lower.7ex\hbox{E}\kern-.125emX}}

\def\BibTeX{{\rm B\kern-.05em{\sc i\kern-.025em b}\kern-.08em
    T\kern-.1667em\lower.7ex\hbox{E}\kern-.125emX}}

\begin{document}

\title{A nonlinear weighted anisotropic total variation regularization for electrical impedance tomography}
\author{Yizhuang Song, Yanying Wang and Dong Liu
\thanks{This work was supported by Shandong Provincial Outstanding Youth Fund (Grand No. ZR2018JL002), National Natural Science Foundation of China (Grand No. 11501336 and No. 61871356) and the China Postdoctoral Science Foundation (2019T120604, 2018M630795). (Corresponding author: Dong Liu.)}
\thanks{Yizhuang Song and Yanying Wang are with School of Mathematics and Statistics, Shandong Normal University, Jinan, Shandong, 250014, P. R. China. Yizhuang Song is also with Center for Post-doctoral Studies of Management Science and Engineering, Shandong Normal University, Jinan, Shandong, 250014, P. R. China (e-mail: ysong@sdnu.edu.cn).}
\thanks{Dong Liu is with CAS Key Laboratory of Microscale Magnetic Resonance and School of Physical Sciences, University of Science and Technology of China, Hefei 230026, China, also with School of Biomedical Engineering and Suzhou Institute for Advanced
Research, University of Science and Technology of China, Suzhou, China, and also with CAS Center for Excellence in Quantum Information and Quantum Physics, University of Science and Technology of China, Hefei 230026, China.  (e-mail: dong.liu@outlook.com).}}



\maketitle

\begin{abstract}
This paper proposes a nonlinear weighted anisotropic total variation (NWATV) regularization technique for electrical impedance tomography (EIT).
The key idea is to incorporate the internal inhomogeneity information (e.g., edges of the detected objects) into the EIT reconstruction process,
aiming to preserve the conductivity profiles (to be detected).
We study the NWATV image reconstruction by employing a novel soft thresholding based reformulation included in the alternating direction method of multipliers (ADMM).
To evaluate the proposed approach, 2D and 3D numerical experiments and human EIT lung imaging are carried out.
It is demonstrated that the properties of the internal inhomogeneity are well preserved and improved with the proposed regularization approach, in comparison to traditional total variation (TV) and recently proposed fidelity embedded regularization approaches.
Owing to the simplicity of the proposed method, the computational cost is significantly decreased compared with the well established primal-dual algorithm.
Meanwhile, it was found that the proposed regularization method is quite robust to the measurement noise, which is one of the main uncertainties in EIT.
\end{abstract}


\section{Introduction}
\label{sec:introduction}
\IEEEPARstart{E}{lectrical} Impedance Tomography (EIT) aims to reconstruct the (change of) conductivity distribution inside objects by injecting a current and measuring the voltage responses through pairs of surface electrodes mounted on the object.
EIT has the advantages of being noninvasive, portable, low cost, capable of high temporal resolution, long duration and continuously monitoring, and much more. These advantages make EIT useful for bedside medical apparatus in clinical applications.
For this, EIT was commercialized and introduced in medical applications since 1980s \cite{Barber1984}.
However, EIT has not yet been widely used in routine clinical applications due to the fact that it is a diffusive modality.

The EIT reconstruction process is commonly recasted into a (least square based) data-fitting inverse problem between the boundary measurement and the computational data.
To deal with the ill-posedness, regularization techniques are widely added to the data-fitting to attract the solution satisfying applications driven constraints. Depending on the form of different constraints, the regularization methods can be roughly classified into three categories: projection based regularization, {\it a priori} conductivity based penalization and learning based regularization.

For projection based regularization, typical examples are truncated singular value decomposition (tSVD) method \cite{Hansen1987,Tehrani2012} and principle component analysis (PCA) method \cite{Vauhkonen1998}.
Even though these methods are capable of providing stable reconstructions, they generally produce ringing artifacts due to the dropout of a certain frequency components \cite{Choi2014}.

In the case of penalty based regularization, typical examples include Tikhonov regularization \cite{Cheney1999,Vauhkonen1998} and its multiplicative form \cite{Zhang2019}, monotonicity based regularization\cite{Zhou2018}, factorization based regularization \cite{Choi2014}, sparsity based regularization \cite{Wang2019a,Wang2020,Shi2020}.
These methods are able to provide stable reconstructions at the cost of blurring the edges of internal inhomogeneities.
Total variation (TV) based regularization \cite{Rudin1992,Borsic2010} and its variants \cite{Li2019,Shi2021a,Shi2021b} have the advantage of preserving the discontinuities of the internal structures especially for dealing with the cases of piecewise constant conductivity distributions.
Since TV regularization is non-differentiable, the so called primal-dual algorithm \cite{Borsic2010} is usually used to deal with the non-differentiability. However, this method needs to handle two optimization problems and hence, it is time consuming \cite{Wu2010} and could produce pseudo-edges and lead staircase effects in the reconstruction. Meanwhile, the Split Bregman method and the first-order TV regularization \cite{Jung2015} are also used at the cost of decreasing the effect of edge preserving \cite{Zhou2015}.
Lee {\it et al} proposed a so-called fidelity embedded regularization (FER) technique \cite{Lee2018}. Using this method, high quality geometries of the internal inhomogeneities can be obtained.
However, since the regularization in the method does not depends on the internal structures the accuracy of values of the reconstructed conductivity, which could be a useful information for clinical use, can not be guaranteed. In actuality, numerical simulations show that when the regularization parameter was set to be infinity, the estimated conductivity is usually far away from the true value (see Section \ref{experiments} for details).
Recently, regularization based on manifold learning \cite{Seo2019} has been published using the results in \cite{Lee2018} as the training data.
Since EIT has not been widely used in clinical applications, it is difficult to obtain enough training data to improve the performance continuously.

In this work, we proposed a nonlinear weighted anisotropic TV regularization approach for EIT reconstruction problem.
In comparison to the well established isotropic TV regularization methods, the nonlinear weighted anisotropic regularization takes use of a nonlinear weight related to the unknown conductivity to eliminate the pseudo edges of the internal structures. Moreover, a soft thresholding formula is derived in the ADMM \cite{Boyd2011} type algorithm to accelerate the reconstruction. Comparing with the anisotropic TV regularization method \cite{Borsic2010,Jung2015}, the proposed regularization approach utilizes nonlinear weight to pull back the internal edges from a possible distortions along the coordinate axes \cite{Gonzalez2017} and provides more accurate EIT reconstructions.
To contextualize the proposed regularization among contemporaries,
Table \ref{ComparsionReg} provides a comparison of the pros and cons
associated with the existing EIT regularization methods in the field of EIT.
To validate the proposed method, we conduct 2D and 3D numerical simulations and human EIT lung imaging experiment.
These experiments are carried out to illustrate the main advantages of the proposed NWATV based method with respect to the existing TV, first-order TV and FER regularization methods.

The remaining sections of this paper are organized as follows. In Section II we review the forward and inverse problems in EIT. Following we describe the proposed nonlinear weighted anisotropic TV regularization approach in Section III and provide the reconstruction algorithm in Section IV. Next we provide the numerical and human lung experiments in Section V. Finally, we provide a discussion in Section VI and conclude the paper in Section VII.

\begin{table*}[]
\centering
\caption{Comparisons of the pros and cons of different regularization schemes in EIT}
\begin{tabular}{p{1cm}p{2.3cm}p{3.5cm} p{4.5cm}p{4.5cm}}
\hline\hline
   Categories & Name of regularizer  & Priors & Pros     & Cons \\ \hline
 \multirow{11}{*}{\rotatebox{90}{\bf Projection method}} & tSVD method \cite{Hansen1987,Tehrani2012}& Singular value of $\mathbf{S}>$a threshold & Easy to implement. & There is ringing artifacts abandoning the basis corresponding to the small singular value \cite{Choi2014}.\\
  &  Principle component analysis (PCA) \cite{Vauhkonen1998} & $\sigma$ lies in the solution space $S_\omega$ & Apply geometrical information from MRI and CT etc. & Need a learning set obtained from the other modalities such as CT, MRI or the reference conductivity value which may not be available. Instability could occur when the priors are not from the patient. \\
   & & & & \\
  \hline
  & & & & \\
  \multirow{48}{*}{\rotatebox{90}{\bf \emph{a priori} conductivity based penalty}}  & Tikhonov \cite{Cheney1999,Vauhkonen1998} & $\| \sigma\|_{L^2}$& The optimization problem is convex and differentiable. & Reconstruction results depends on an artificially chosen regularization parameter. Inclusions is blurred. \\
  & Multiplicative Tikholnov regularization \cite{Zhang2019} & $\left\| \frac{|\na\sigma|^2+\delta_{n-1}^2}{|\na\sigma^{n-1}|^2+\delta_{n-1}^2}\right\|_{L^2}$  for measurement data and mesh dependent $\delta_{n-1}$& Reconstruction results are independent of the regularization parameter choice. & A weighted $L^2$ regularizer could still blur the image. In each step, numerical integral and differentiation are needed to calculate. Moreover, Gauss-Newton method with step size using line search method is used to solve a nonquadratic cost functional which is time consuming. \\
  & Total variational (TV) regularization \cite{Borsic2010} & $\int_\Om |\na\sigma|$ & Can preserve the edge of internal inhomogeneities & Staircase effect and pseudo edge exists.  Reconstructing is slow and resolution is high when the prime-dual algorithm is used while resolution is low and reconstructions is fast when the Split-Bregman method is used to minimized the TV constrained optimal problem \cite{Zhou2015}.\\
  & TGV \cite{Shi2020,Gong2018} & $TGV_G^\alpha(\sigma)$ & Efficiency in removing the staircase effect when using TV for piecewise linear conductivity distribution. & In some imaging cases such as human lung imaging, piecewise linear is not a reasonable assumption. Balancing efficiency and quality is difficult as that in TV regularization.   \\
  & First order TV regularization \cite{Jung2015,Gonzalez2017} &  $\|\na\sigma\|_{L^1(\Om)}$ & The imaging speed is high & The capacity of edge preserving is lower than TV. Distort the inclusion boundaries along the coordinate axes. \\
  & Monotonicity based method \cite{Zhou2018} & The sign of $\frac{\partial \sigma}{\partial t}$ coincides with the breathing process. & No confusion between breathing and expirtion & Still a kind of $L^2$ regularization which could blur the internal edges.\\
  & Factorization based method \cite{Choi2014} & Highlight the right singular vector associate with high singular values. &Alleviated Gibbs ringing artifacts & Heuristical argument without strict mathematical theory\\
    & Wavelet frame TGV \cite{Shi2021b}& $TGV+ \|{\bm \lambda}\cdot W \sigma_{recon}\|_0$  & Using the $L^0$ norm of the wavelet frame to sharpen the TGV reconstructed conductivity & Primal-dual method is used for the TGV minimization, reconstruction is time consuming. Postprocessing the TGV reconstructions needs more time.
   \\
   & Fidelity-imbedded regularization \cite{Lee2018} & $\sqrt{\sum\limits_{l}|\langle S_k,S_l \rangle|}$ & Reconstruction results are independent of the regularization parameter choice and fast imaging due to only a direct algebraic inversion is needed to compute.& The accuracy of the reconstruction can not be guaranteed due to the lack of mathematical background. \\
   & & & & \\
   \hline
   & & & & \\
  \multirow{1}{*}{\rotatebox{90}{\bf Learning method}} &Manifold learning method \cite{Seo2019} & $\sigma\in $ a manifold $\mathcal{M}$. & Imaging quality is higher than the model based regularizer. & To obtain experimental labelled data $(\{\mathbf{V}_n,\sigma_n)\}$ is difficult before designing a high quality model based regularization method. \\
   & & & & \\
   & & & & \\
   & & & & \\
   \hline
  & & & & \\
   \multirow{1}{*}{\rotatebox{90}} & {\bf Proposed method} & $\left\|\f{\na\sigma}{|\na\sigma|^2}\right\|_{L^1}$ & Improved preservation of the internal inhomogeneity edges; the computational time is significantly decreased; robust to the additional Gaussian random noise. & Need to select proper regularization parameters manually. \\
& & & & \\
    \hline\hline
\end{tabular}
\label{ComparsionReg}\\
\end{table*}

\section{Forward and inverse problems in EIT}

Let $\Omega \subset \mathbf{R}^n$ ($n=2,3$) be the imaged object with a smooth boundary. Let $\gamma = \sigma+i\omega\epsilon$ represent the admittivity distribution of the region $\Omega$. Here, the conductivity $\sigma$ and permittivity $\epsilon$ depend on position $\r = (x,y,z)$, angular frequency $\omega$.

For an $E$-channel EIT system, $E$ electrodes $\mathcal{E}_1, \mathcal{E}_2, \cdots, \mathcal{E}_E$ are attached on $\partial \Omega$, the boundary of $\Omega$. We inject a series of time harmonic currents with magnitude $I$ mA following e.g. the neighboring protocol. Under such a protocol, we sequentially inject several currents using pairs of electrodes $(\mE_j,\mE_{j+1})$ for $j=1,...,E$, where we assign $E+1$ to be $1$.
The induced electric potential $u^{j}$ is governed by the following elliptic partial differential equation (PDE) with mixed boundary conditions:
\begin{equation}\label{forward}
  \left\{
    \begin{split}
       & \nabla\cdot(\sigma\nabla u^{j})=0,\quad \mbox{ in }  \Omega \\
       & \sigma\nabla u^j\cdot \mathbf{n}=0, \quad \mbox{ on }  {\partial \Om \setminus \mathop{\cup}\limits_{i=1}^{E}\mE_i} \\
       & \int_{\mathcal{E}_j} \sigma \frac{\partial u^j}{\partial {\bf n}}dS = I = -\int_{\mathcal{E}_{j+1}} \sigma \frac{\partial u^{j}}{\partial {\bf n}}dS \\
       & \left. u^j + z_i \sigma\frac{\partial u_j}{\partial {\bf n}}\right|_{\mE_i} = U^j_i,\qquad i, j=1,2,\cdots,E \\
       & \int_{\mathcal{E}_k} \sigma \frac{\partial u^j}{\partial {\bf n}}dS = 0 ,\quad k\in \{1,2,\cdots,E\}\setminus \{j,j+1\}.
    \end{split}
  \right.
\end{equation}
Here, $z_i$ is the contact impedance between the electrode and $\partial \Omega$ and $U_i^j$ is the voltage potential on $\mE_i$ caused by the injection currents using the electrode pair $(\mE_j,\mE_{j+1})$. Using an EIT measurement device, the quantity $V_i^j[\sigma] = U_{i}^j[\sigma] - U_{i+1}^j[\sigma]$ is measurable.
Note that $z_i$ is unknown, however, for $i\in \{1,2,\cdots,E\}\setminus \{j-1,j,j+1\}$, the contact impedance can be neglected \cite{Lee2018} and hence $V_i^j \approx u^j|_{\mE_i} - u^j|_{\mE_{i+1}}$. To reconstruct $\sigma$, EIT uses the following reciprocity principle
\begin{equation}\label{reciprocity}
  V_i^j[\sigma] \approx u^j|_{\mE_i} - u^j|_{\mE_{i+1}} =\frac{1}{I}\int_\Om \sigma(\r) \nabla u^i(\r) \cdot \nabla u^j(\r)d\r.
\end{equation}
EIT uses $E\times (E-3)$ data $V_i^j$ and the relation \eqref{reciprocity} to reconstruct the unknown quantity $\sigma$. However, the above problem is nonlinear and ill-posed. Linearization is usually applied to deal with the non-linearity and regularization is used to handle the ill-posedness.

To be precise, we assume that $\sigma$ is a perturbation of a constant $\sigma_0$, that is $\sigma(\r) = \sigma_0 + \delta \sigma(\r)$ for $\r\in \Omega$. Here, $\delta\sigma$ has compact support in $\Omega$. Then $\delta\sigma$ approximately satisfies the following integral equation
\begin{equation}\label{reci_diff}
  \int_\Omega \delta\sigma(\r)\nabla u^i_0(\r)\cdot \nabla u^j_0(\r)d\r = I(V_i^j[\sigma]-V_i^j[\sigma_0]) =:I\delta V_i^j[\sigma],
\end{equation}
where $u^i_0$ is the solution of the equation \eqref{forward} with $\sigma$ replaced by $\sigma_0$.

Discretize the domain $\Om$ as $\Omega = \cup_{k=1}^N T_k$, where $T_k$ is triangular element, then $\sigma|_{T_k}$ can be considered as a constant $\sigma_k$. $\sigma$ can be approximated by the vector $\bm{\sigma}:=(\sigma_1,\sigma_2,\cdots,\sigma_N)^T$. Then we convert the formula \eqref{reci_diff} to the following system
\begin{equation}\label{eq:sensitivity}
\mathbf{S}{\delta \bm \sigma} = \delta \mathbf{ V}
\end{equation}
where $\mathbf{S}$ is an $E(E-3)\times N$ \emph{sensitivity matrix} whose element $S_{pq} = \frac{1}{I} \int_{T_q} \na u_0^m(\r) \cdot \na u_0^n(\r)d\r$ , where $m = [\frac{p-1}{N-3}]+1$, $n = {p-12[\frac{p-1}{N-3}]+2}$,
$\delta \mathbf{V}\in \mathbf{R}^{E(E-3)}$ whose $q$-th element is $V_q = \delta V_m^n$
and the operator $[\cdot]$ represents the largest integer less than $\cdot$.

\section{Nonlinear weighted anisotropic TV regularization}

To solve the equation \eqref{eq:sensitivity}, we reformulate it to the following least squares problem
\begin{equation}\label{least-square}
  \delta\bm{\sigma}^* = \mathop{{\rm arg\ min}}\limits_{\delta\bm{\sigma}}\|\mathbf{S}\delta\bm{\sigma} - \delta\mathbf{V} \|^2_{L^2(\Omega)}.
\end{equation}

Let $\mathbf{s}_l$ denotes the $l$-th column vector of $\mathbf{S}$ for $l=1,2,\cdots,N$. The correlation between the $l$ and $k$-the column vectors $\langle \frac{\mathbf{s}_l}{|\mathbf{s}_l|}, \frac{\mathbf{s}_k}{|\mathbf{s}_k|}\rangle$ decreases rapidly as the distance between $T_l$ and $T_k$ increases \cite{Lee2018}; this makes the condition number of the matrix $\mathbf{S}^T \mathbf{S}$ to be approximately infinity. Hence, the minimization problem \eqref{least-square} is unstable against the measurement errors and noises in $\delta \mathbf{V}$.
We use the following minimization problem to approximate \eqref{least-square},
\begin{equation}\label{regularization}
  \delta\bm{\sigma}^*_{\lambda} = \mathop{{\rm arg\ min}}\limits_
{\delta\bm{\sigma}}\left\{\frac{1}{2}\|\mathbf{S}\delta\bm{\sigma} - \delta\mathbf{V} \|^2_{L^2(\Omega)}+ \lambda\mbox{Reg}(\delta\bm{\sigma})\right\}.
\end{equation}
Here, instead of minimizing the term $\|\mathbf{S}\delta\bm{\sigma} - \delta\mathbf{V} \|^2_{L^2(\Omega)}$, we also minimize the other regularization term $\mbox{Reg}[\delta{\bm{\sigma}}]$, where $\lambda$ is a parameter balancing the fidelity term $\|\mathbf{S}\delta\bm{\sigma} - \delta\mathbf{V} \|^2_{L^2(\Omega)}$ and the regularization term. The choice of regularization depends on the {\it a priori} information of the conductivity.

Following, for ease of explanation we write $\Omega = \cup_{i=1}^{n_I} D_i \cup (\Omega\setminus \overline{\cup_{i=1}^{n_I} D_i})$, where $D_i$ ($i=1,2,\cdots,n_I$) represent $n_I$ internal structures with $n_I$ a positive integer.
Notice that $\zeta(\r) := 1/|\nabla \sigma(\r)|^2$ can be considered as indicators of edges of conductivity inhomogeneities. To be precise, as ${\bf r}\to \partial D_i$, $\zeta(\r)\to 0$ and $\zeta(\r)\approx +\infty$ for $\r$ away from $\partial D_i$ ($i=1,2,\cdots,n_I$).
 Moreover, $|\nabla \sigma|$ is relatively sparse in lung imaging due to the piecewise constant structure of $\sigma$ at a given time.
Based on the above observations, we construct the nonlinear weighted regularization term $\|\nabla\sigma/|\nabla\sigma|^2\|_{L^1(\Omega)}$. Hence, \eqref{regularization} is changed to
\begin{equation}\label{proposed method}
 \delta\bm{\sigma}^{*}_{\lambda}=\mathop{{\rm arg\ min}}\limits_{\delta\bm{\sigma}}  \left\{\frac{1}{2}\|\mathbf{S}\delta\bm{\sigma}-\delta\bm{V}\|^{2}+\lambda\|\mathbf{p}\cdot\mathbf{D}(\delta\bm{\sigma})\|_{l_1}\right\},
\end{equation}
Here, $\mathbf{p}=(\zeta(\delta\bm{\sigma});\zeta(\delta\bm{\sigma}))\in \mathbf{R}^{2N}$, $\mathbf{D}(\delta\bm{\sigma}) = (\mathbf{D}_{x}{\delta\bm{\sigma}};\mathbf{D}_{y}{\delta\bm{\sigma}})\in \mathbf{R}^{2N}$,
where $\mathbf{D}_{x}$, $\mathbf{D}_{y} \in \mathbf{R}^{N\times N}$ are respectively the first order difference operators along the $x$ and $y$ directions.

Given the fact that the minimization problem \eqref{proposed method} is non-differentiable, we adopt the well established Alternating Direction Method of Multipliers (ADMM) \cite{Boyd2011} to solve it. The initial guesses $(\mathbf{z}_0,\mathbf{p}_0,\mathbf{y}_0)$, $\delta\bm{\sigma}$ are updated via the following schemes,
\begin{subequations}
\begin{align}
\label{a5}
\delta\bm{\sigma}_{n+1}&=\mathop{{\rm arg\ min}}\limits_{\delta\bm{\sigma}}\mathcal{L}_{\rho}(\delta\bm{\sigma}, \mathbf{z}_{n}, \mathbf{p}_{n}; \mathbf{y}_{n});\\
\label{a6}
\mathbf{z}_{n+1}&=\mathop{{\rm arg\ min}}\limits_{\mathbf{z}}\mathcal{L}_{\rho}(\delta\bm{\sigma}_{n+1}, \mathbf{z}, \mathbf{p}_{n}; \mathbf{y}_{n});\\
\label{a7}
 \mathbf{p}_{n+1}&=(\zeta(\delta\bm{\sigma}_{n+1});\zeta(\delta\bm{\sigma}_{n+1}));\\
\label{a8}
\mathbf{y}_{n+1}&=\mathbf{y}_{n}+\rho(\mathbf{D}\delta\bm{\sigma}_{n+1}-\mathbf{z}_{n+1}).
\end{align}
\end{subequations}
Here, $\mathbf{z}\in \mathbf{R}^{2N}$ is an auxiliary variable, $\mathcal{L}_\rho$ is the augmented Lagrangian functional defined as
\begin{equation*}\label{Lagrangian functional 1}
\begin{split}
\mathcal{L}_{\rho}(\delta\bm{\sigma}, \mathbf{z}, \mathbf{p}; \mathbf{y}):=\frac{1}{2}\|\mathbf{S}\delta\bm{\sigma}-\delta \mathbf{V}\|^{2}+\lambda\|\mathbf{p}\cdot \mathbf{z}\|_{{l_1}}\\
+\frac{\rho}{2}\|\mathbf{D}\delta\bm{\sigma}-\mathbf{z}\|_{l_2}^{2}+\mathbf{y}^{T}(\mathbf{D}\delta\bm{\sigma} - \mathbf{z}),
\end{split}
\end{equation*}
where $\mathbf{y}$ is the Lagrangian multiplier and $\rho$ is a scalar penalty parameter.

Meanwhile, from \eqref{a5} and \eqref{a6} we can derive
\begin{equation}
\label{a51}
\delta\bm{\sigma}_{n+1}=(\frac{1}{\rho}\mathbf{S}^{T}\mathbf{S}+\mathbf{D}^{T}\mathbf{D})^{-1}(\frac{1}{\rho}\mathbf{S}^{T}\delta\mathbf{V}+\mathbf{D}^{T}\mathbf{z}_{n}-\mathbf{D}^{T}\frac{\mathbf{y}_{n}}{\rho}),
 \end{equation}
and
\begin{equation}\label{a62}
\mathbf{z}_{n+1}[k]=h_{\frac{\lambda\mathbf{p}_{n}[k]}{\rho}}(\mathbf{D}\delta\bm{\sigma}_{n+1}[k]+\frac{\mathbf{y}_{n}[k]}{\rho}).
\end{equation}
Here, $\mathbf{z}_{n+1}[k]$ is the $kth$ element of $\mathbf{z}_{n+1}$ for $k=1,2,\cdots,2N$, $h_{g}(\cdot)$ is the soft thresholding formula defined as \cite{Daubechies2004}
 \begin{equation}\label{def:h_g}
 h_{g}(\cdot)=
  \left\{\begin{array}{ll}
  \cdot-g\mbox{sgn}(\cdot),&|\cdot|>g,\\
0, & \mbox{otherwise},
\end{array}\right.
\end{equation}
where $\mbox{sgn}$ is the sign function. We provide a proof of formula \eqref{a62} in the Appendix section.

\section{Reconstruction algorithm}

In this section, we summarize the reconstruction algorithm utilizing the proposed nonlinear weighted anisotropic TV regularization. Note that in the human experiment there is unavoidable noise and artifacts in the measured data, we first pre-process the data as that in \cite{Lee2018} to reduce the artifacts caused by rib cage movement\cite{liu2015nonlinear}.

\subsection{Data pre-processing and modeling error blocking}
Instead of using the direct measurement $\delta\mathbf{V}$, we use the pre-processed data $\delta\mathbf{V}^\diamond$.
Here,
\begin{equation*}
\delta \mathbf{V}^\diamond = \delta\mathbf{V} - \mathbf{S}_{bdy}(\mathbf{S}^{T}_{bdy}\mathbf{S}_{bdy}+\lambda_{b}\mathbf{I})^{-1}\mathbf{S}_{bdy}^{T}\delta \mathbf{V},
\end{equation*}
where $\mathbf{S}_{bdy}$ is a sub-matrix of $\mathbf{S}$ including the columns related with boundary elements, $\lambda_b$ is a regularization parameter and $\mathbf{I}$ is the identity matrix.

Due to the reconstruction problem is iteratively solved, the modeling error could propagate in the forward problem, especially the estimated conductivity near the electrodes during the iterations.
The modeling error propagation could heavily deteriorate the accuracy of the reconstructed images.
To this end, we assume that the conductivity near the electrodes is invariant with respect to patients breathing and it only changes in the domain $\widetilde \Om$, where $\widetilde \Om \subset\subset \Om$ represents the domain including two lungs. For each iteration, we modify $\delta{\bm\sigma}_{n+1}$ by
\begin{equation}\label{modify_result}
\delta {\bm \sigma}_{n+1} \leftarrow \delta {\bm \sigma}_{n+1}\chi_{\widetilde \Om},
\end{equation}
where $\chi_{\widetilde \Om}$ is the characteristic function of the domain $\widetilde \Om$. Subsequently, we will nominate modified nonlinear weighted anisotropic TV (MNWATV) as the proposed regularization using the preprocessed data.

\subsection{Peudocode for the proposed method}
We summarize the above procedures as the following algorithm in the form of pseudocode.

\begin{algorithm}[h]
  \caption{The proposed algorithm }
  \label{The proposed algorithm}
  \begin{algorithmic}[1]
    \REQUIRE {Measured data $\mathbf{V}$, the maximum iteration number $M$, the parameters $\lambda$, $\rho$, tolerance $tol$ and the sensitivity matrix $\mathbf{S}$}
    \ENSURE {conductivity distribution $\bm{\sigma}$}
    \STATE {Initialize: $\mathbf{y}_{0}=\mathbf{0}$, $\mathbf{z}_{0}=\mathbf{0}$, $\mathbf{p}_{0}=\mathbf{1}$, $n=1$.}
    \WHILE {$n<M$}
    \STATE {Update $\bm{\sigma}_{n+1}$ using \eqref{a51} and set $\bm{\sigma}_{n+1} \leftarrow \bm{\sigma}_{n+1}\chi_{\widetilde \Om}$}.
    \STATE {Update $\mathbf{z_{n+1}}$ using \eqref{a62}.}
    \STATE {Update $\mathbf{p}_{n+1}$ using \eqref{a7}.}
    \STATE {Update $\mathbf{y}_{n+1}$ using \eqref{a8}.}
    \IF{$\|\bm{\sigma}_{n+1}-\bm{\sigma}_{n}\|<tol$} \STATE{break} \ENDIF
    \ENDWHILE
  \end{algorithmic}
\end{algorithm}

\section{Numerical and Human Experiments}\label{experiments}
To test the performance of the proposed algorithm, we carried out numerical examples and human EIT lung imaging experiment.
\subsection{Experiments Setup}
In the experiments, we attach 16 electrodes to the boundary of the subject as uniform as possible. In all the experiments, we inject currents with an amplitude of 1.0 mA and measure the voltages between adjacent pairs of electrodes. To avoid having $0$ in the denominator, in all the experiments we approximate $\zeta(\delta \bm\sigma_{n+1})$ by $\f{1}{|\D \delta\bmsigma_{n+1}|^2+\delta}$ for a small $\delta>0$.

The results of the experiments are compared to methods employing TV, first-order TV and FER regularizers. For the comparison, in the numerical experiments, we computed the relative error (RE), the peak signal-to-noise ratio (PSNR) to evaluate the accuracy of reconstructed images.
At the $n$-th step, $RE(n)$ is defined as
\begin{equation*}\label{RE}
RE(n)=\frac{\|\bm{\sigma}_n-\bm{\sigma}^{\ast}\|_2}{\|\bm{\sigma}^{\ast}\|_2},
\end{equation*}
and $PSNR(n)$ is defined as
\begin{equation*}
\label{psnr}
PSNR(n) = 10 \log_{10}\frac{\max(\bm{\sigma}_{n}\odot \bm{\sigma}_{n})}{MSE(n)}.
\end{equation*}
Here $\cdot\odot\cdot $ represents the component wise multiplication. $\bm{\sigma}^{\ast}$ represents the true conductivity and $MSE(n) = \frac{1}{N^2}\sum_{i=1}^{N^2}\left[\bm{\sigma}_n(i)-\bm{\sigma}^\ast(i)\right]^{1/2}$ is the mean square error between $\bm{\sigma}^{\ast}$ and $\bm{\sigma}_{n}$.

Since we do not know the true conductivity in the human lung imaging experiment,
to define the $\widetilde{RE}(n)$ and $\widetilde{PSNR}(n)$,
$\sigma^*$ was replaced by $\sigma_{TV}$, which is the reconstructed conductivity using the TV regularization.
Additionally, we also compare the CPU time in the reconstruction process.
The reconstructions are carried out using  a workstation with 2.00 GHz Inter (R) Xeon (R) Gold 6138 CPU, 256 GB memory, Windows 10 operating system.

To decrease the computational cost, we use the point electrodes at the center of the true electrodes instead of using the true circular electrodes.
\cite{Hanke2011} provide a theoretical result towards the reasonability of such an approximation. We present the reconstruction results in the image with $256\times 256$ pixels. Moreover, in all the experiments, we set $tol = 10^{-5}$.

\subsection{2D numerical experiments}

\begin{figure}[htbp]
 \centering
  \subfloat[]{\label{model_6}
   \includegraphics[width=0.25\linewidth]{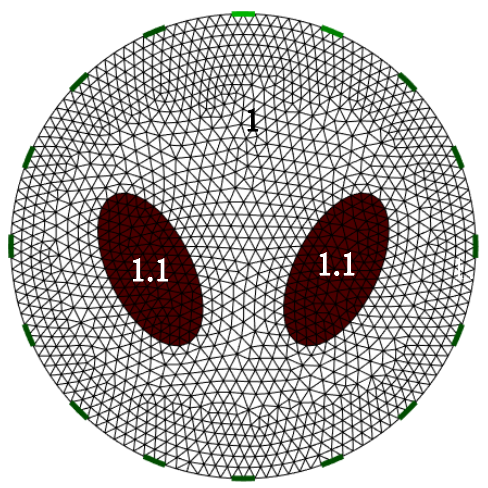}}
   ~~~~~~
  \subfloat[]{\label{Inv_model_2D}
   \includegraphics[width=0.25\linewidth]{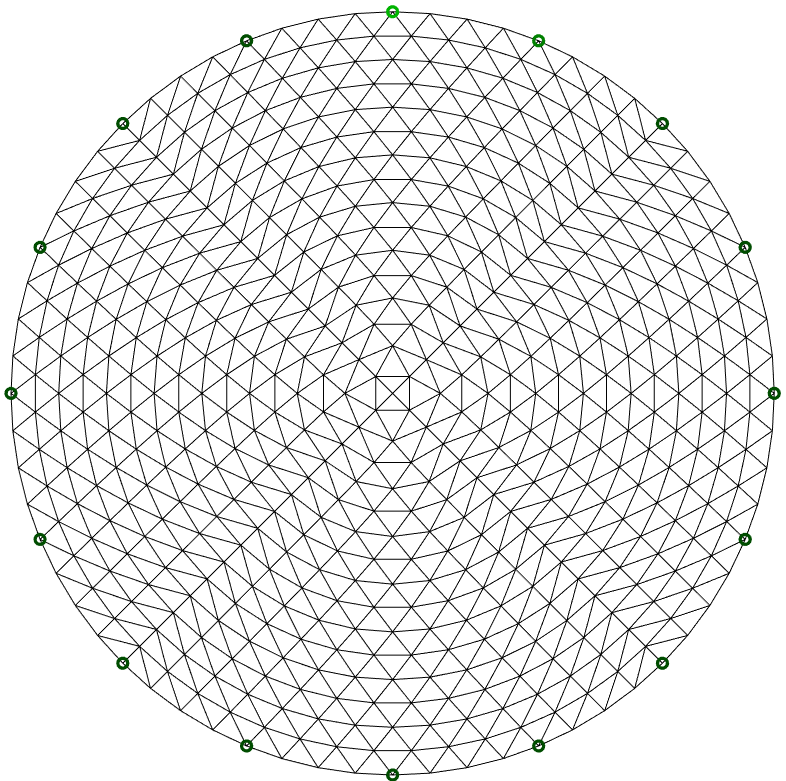}}

 \caption{(a) The 2D numerical model and its conductivity distribution.  (b) Nodes and elements used in the inversion process.}
\label{Numerical_mode_2D}
\end{figure}

\begin{figure*}[htbp]
  \centering
  \includegraphics[width=0.95\linewidth]{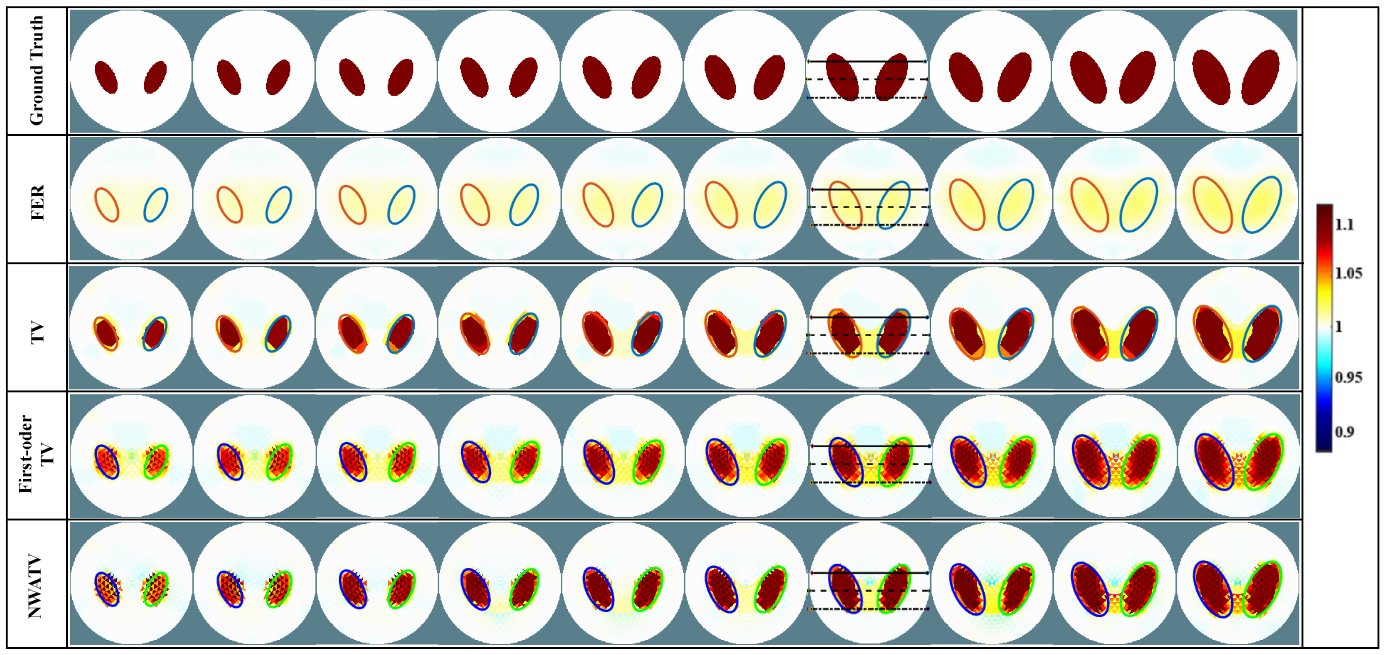}\\
  \caption{The reconstruction results for 2D model using four different regularizers. The first row depicts the ground truth conductivity distribution, the second to the fifth row illustrate the reconstruction results using the FER, the TV, the first-order TV and the proposed regularizers. }\label{lung_results}
\end{figure*}

We construct a disk with radius of $0.1$ m centered at the original point $(0,0)$. We put two ellipses inside the disk to simulate the lungs. The center of the constructed ellipses are respectively $(0.04,-0.01)$ and $(-0.04,-0.01)$. The size of the ellipses vary to form 10 different models (see the first row of Fig. \ref{lung_results}). For the $k$-th model ($k=1,2,\cdots,10$), the major and minor semi-axes of the left ellipse are respectively $\mathbf{a}_l=(0.012+0.001k, 0.024+0.002k)$, $\mathbf{b}_l=(-0.012-0.001k, 0.006+0.0005k)$ while the values of the right ellipse are respectively $\mathbf{a}_r=(-0.012-0.001k, 0.024+0.002k)$, $\mathbf{b}_r=(0.012+0.001k, 0.006+0.0005k)$. We set the conductivity of background and inclusion to be 1.0 S/m and 1.1 S/m respectively. One of the 10 models is shown in Fig. \ref{Numerical_mode_2D} (a). In the inversion, a finite element mesh with 1024 triangular elements and 545 nodes as shown in Fig. \ref{Numerical_mode_2D} (b) is used.

We first solve the two-dimensional forward problem \eqref{forward} to obtain the boundary voltage datum. Using this datum we reconstruct the conductivity images using the proposed algorithm and three existing methods. Gaussian random noise with SNR 50dB is added to the data to test the anti-noise performance.  We set $\lambda=5\times10^{-13},\rho=1\times 10^{-10}, \delta=0.01, M=20$ for the proposed method. We explain the selected parameters in Section \ref{section:discussion}. For the FER method, we only consider the case when the regularization parameter is set to be $\infty$. The reference conductivity $\sigma_0$ is set to be background conductivity 1.0. The parameters of the other regularization methods are empirically chosen.
The results are shown in Fig. \ref{lung_results}.

\begin{figure*}[htbp]
\centering
\subfloat[]{
\includegraphics[width=0.23\linewidth]{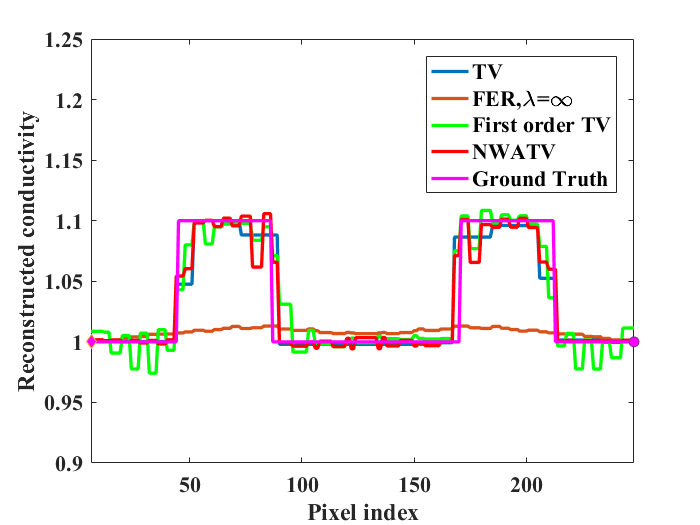}}
\subfloat[]{
\includegraphics[width=0.23\linewidth]{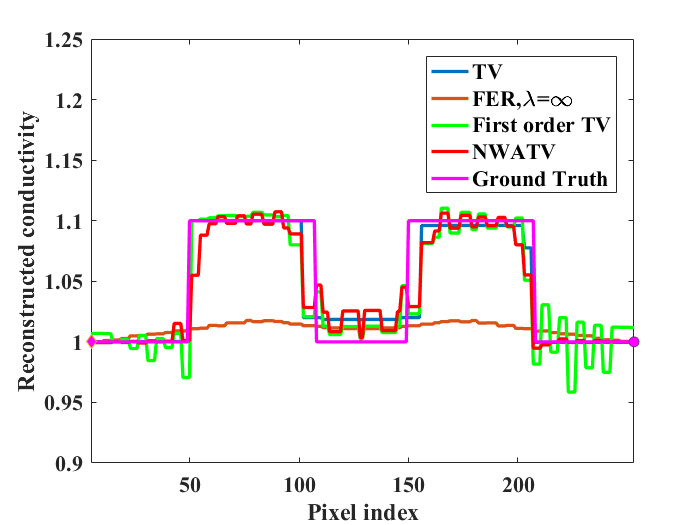}}
\subfloat[]{
\includegraphics[width=0.23\linewidth]{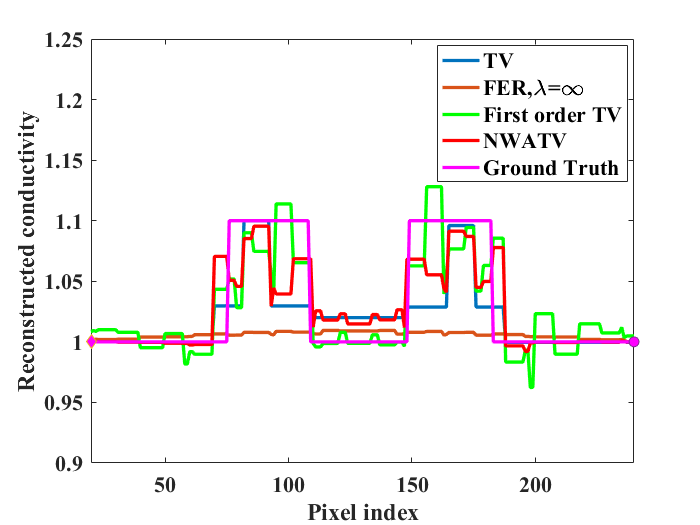}}
\centering
\caption{Profiles of the reconstructed conductivity along the solid line (a), the slash line (b) and the dot dash line (c) shown in Fig. \ref{lung_results}. The meaning of each colorful lines in (a)-(c) are shown in each subfigure.}
\label{1D profiles_lung}
\end{figure*}

\begin{figure}[h]
\centering
\subfloat[]{
\includegraphics[width=0.5\linewidth]{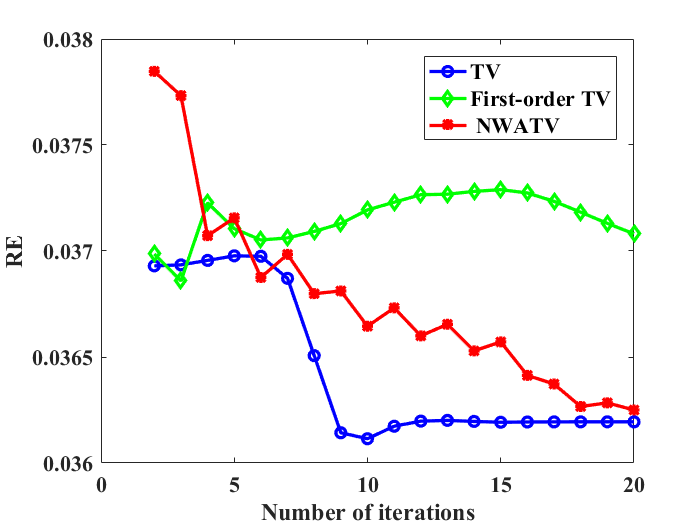}}
\subfloat[]{
\includegraphics[width=0.5\linewidth]{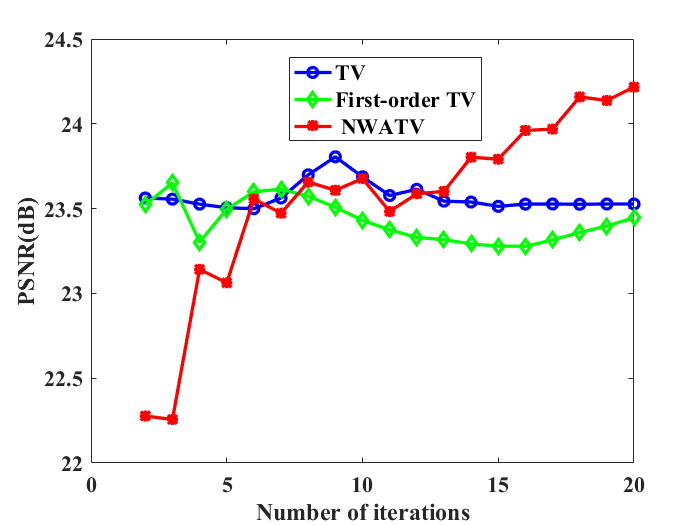}}
\centering
\caption{Behaviors of $RE(n)$ (a) and $PSNR(n)$ (b) of 2D simulation. Meanings of each colorful lines are shown in the subfigure.}
\label{quantitative indicators_2d_lung}
\end{figure}

Fig. \ref{1D profiles_lung} (a)-(c) illustrates the profiles along the solid, the dash and the dot dash lines respectively shown in the 7th model in Fig. \ref{lung_results}. Fig. \ref{quantitative indicators_2d_lung} (a)-(b) shows the behavior of $RE(n)$ and $PSNR(n)$ as $n$ increases for each regularizers except FER since it is a direct reconstruction method. In Table II we compare the computational time for the reconstructions using the aforementioned regularizers.

From the comparisons, the other three regularizers can produce more accurate images than FER does. Moreover, for the proposed method, the computational time is significantly reduced while maintaining a satisfactory accuracy. On the other hand the accuracy of the first-order TV regularizer is obviously decreased.

\subsection{3D numerical experiment}

\begin{figure}[htbp]
 \centering
  \subfloat[]{\label{Side view}
   \includegraphics[width=0.28\linewidth]{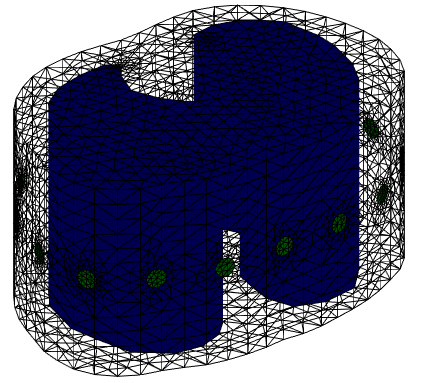}}
 \hfill
  \subfloat[]{\label{conductivity distribution}
   \includegraphics[width=0.28\linewidth]{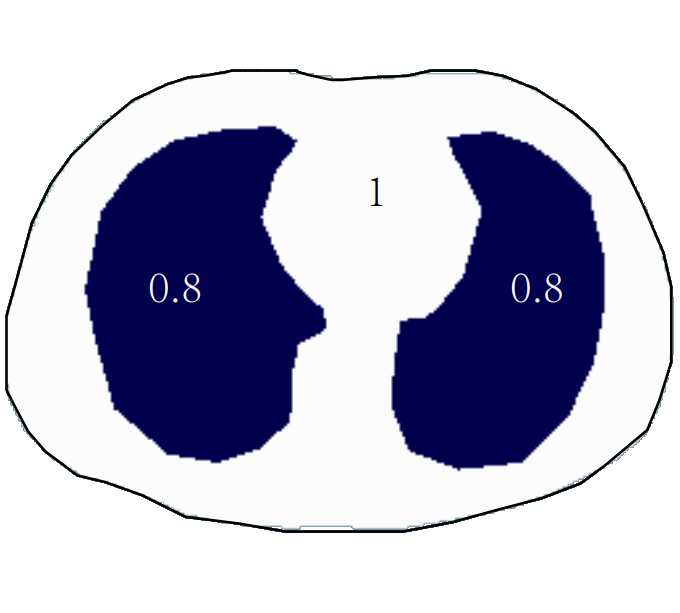}}
 \hfill
  \subfloat[]{\label{Imaging plane}
   \includegraphics[width=0.28\linewidth]{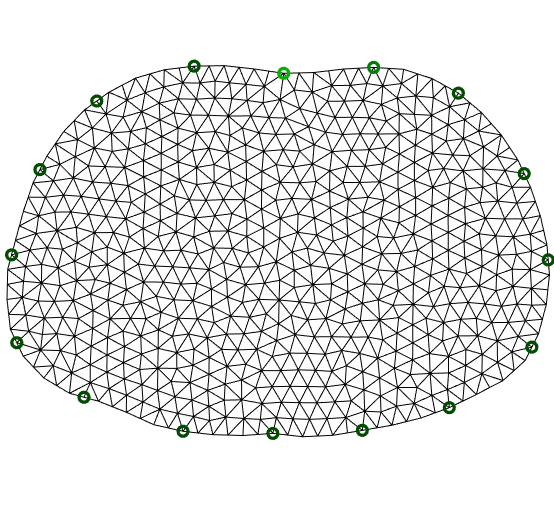}}

 \caption{The 3D adult human thorax model. \protect\subref{Side view} The geometry and internal structure. \protect\subref{conductivity distribution} Ground truth conductivity distribution at the center slice. \protect\subref{Imaging plane} Imaging plane.}
\label{Numerical_model}
\end{figure}

We do the numerical simulation using a 3D adult human thorax model in EIDORS \cite{Adler2006}. The geometry of the model is shown in Fig. \ref{Numerical_model} (a).
We attach 16 circular electrodes along the boundary of the center slice. The geometry is discretized by 5,379 nodes 26,358 tetrahedrons. To test the algorithm, we set a conductivity contrast as follows: the conductivity of the background is set to be 1.0 S/m while the values in the two lung shape anomalies are set to be 0.8 S/m. Fig. \ref{Numerical_model} (b) illustrates the true conductivity distribution at the center slice. Fig. \ref{Numerical_model} (c) provide the discretization information of the reconstruction which includes 637 nodes and 1,172  triangles.

\begin{figure*}[htbp]
\centering
 \subfloat[]{
  \includegraphics[width=0.2\linewidth]{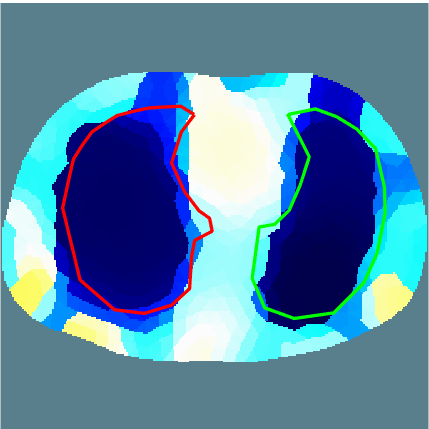}}
 \subfloat[]{
  \includegraphics[width=0.2\linewidth]{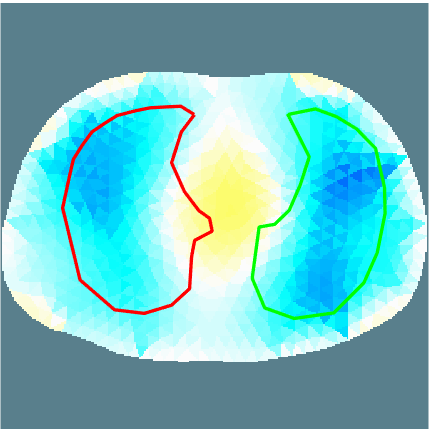}}
 \subfloat[]{
  \includegraphics[width=0.2\linewidth]{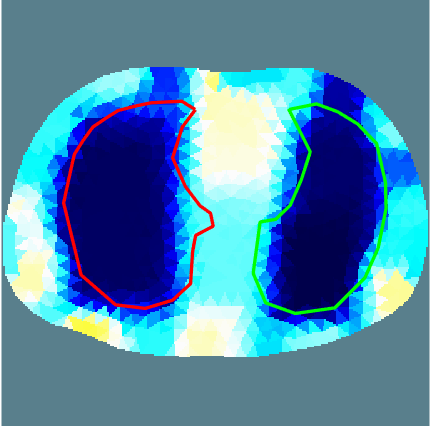}}
 \subfloat[]{
  \includegraphics[width=0.238\linewidth]{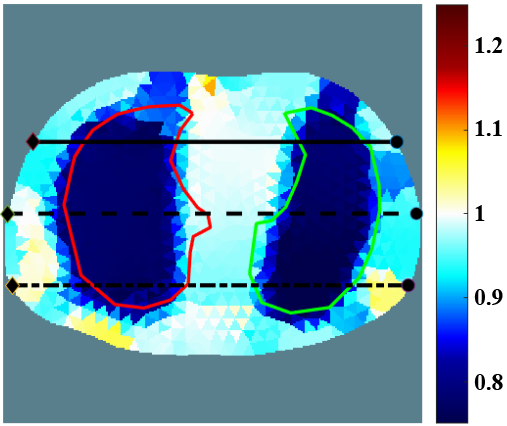}}
  \centering
 \caption{The reconstruction results using four different regularizers. (a) is the reconstruction using the primal-dual algorithm for TV regularizer at the 10th steps. (b) depicts the result of FER for the regularization parameter set to be $\infty$. (c) is the reconstructions using the first-order TV regularizer at the 10th iteration. (d)is the results using the NWATV at the 10th step.}
\label{numerical_results}
\end{figure*}

\begin{figure*}[htbp]
\centering
\subfloat[]{
\includegraphics[width=0.25\linewidth]{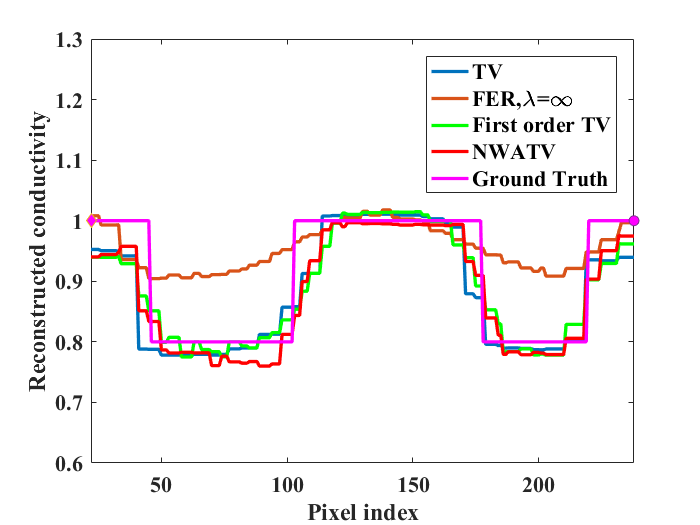}}
\subfloat[]{
\includegraphics[width=0.25\linewidth]{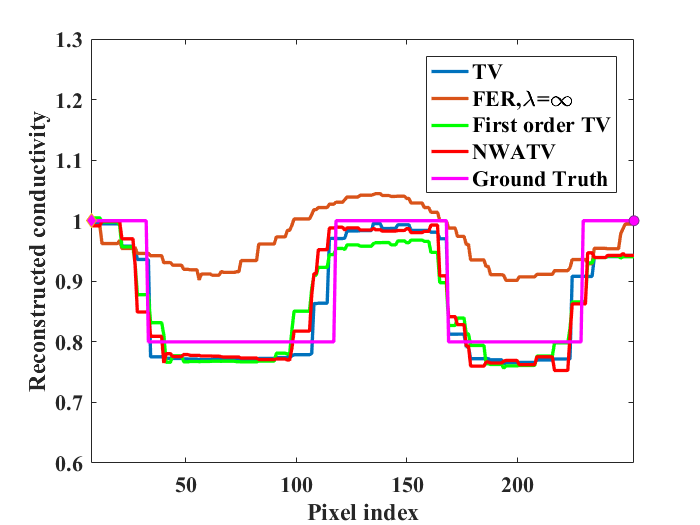}}
\subfloat[]{
\includegraphics[width=0.25\linewidth]{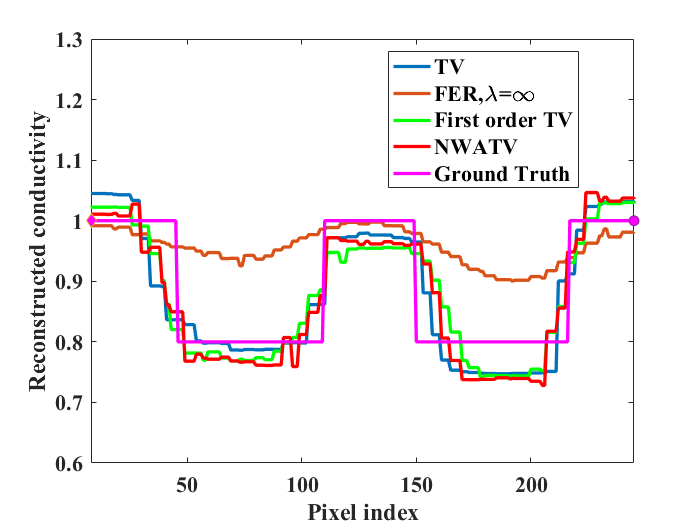}}
\centering
\caption{Profiles of the reconstructed conductivity along the solid line (a), the slash line (b) and the dot dash line (c) shown in Fig. \ref{numerical_results} (d). The meaning of each colorful lines in (a)-(c) are shown in each subfigure.}
\label{1D profiles}
\end{figure*}

\begin{figure}[h]
\centering
\subfloat[]{
\includegraphics[width=0.5\linewidth]{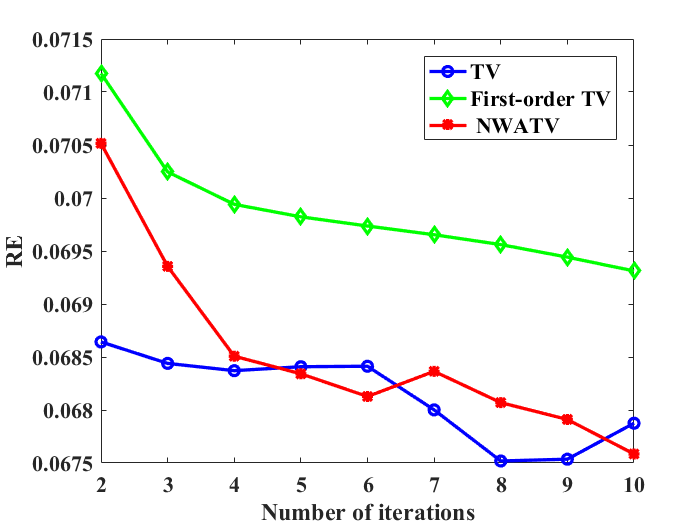}}
\subfloat[]{
\includegraphics[width=0.5\linewidth]{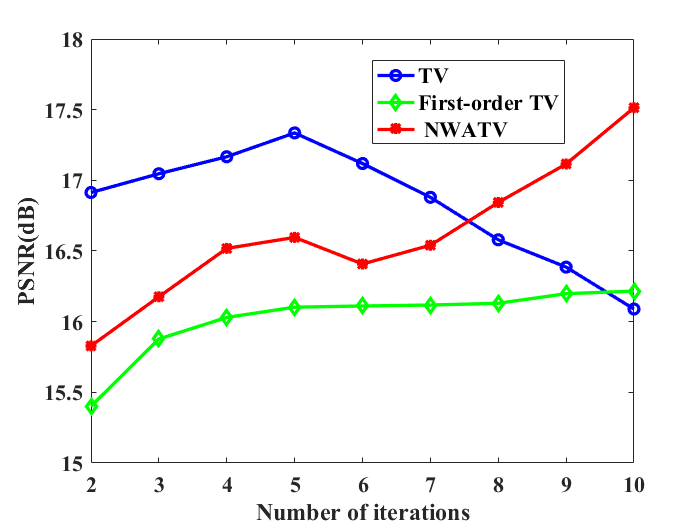}}
\centering
\caption{Behaviors of $RE(n)$ (a) and $PSNR(n)$ for the 3D simulation. Meanings of each colorful lines are shown in the subfigure.}
\label{quantitative indicators}
\end{figure}

We first solve the three-dimensional forward problem \eqref{forward} to obtain the boundary voltage datum. Using this datum we reconstruct the conductivity images using the proposed algorithm and three existing methods. The results are shown in Fig. \ref{numerical_results}. In Fig. \ref{numerical_results}, (a) shows the reconstruction using TV regularization and primal-dual algorithm at the 10th step, (b) depicts the results of FER for the regularization parameter set to be $\infty$, (c) illustrates the reconstructions using the first-order TV regularizer and the ADMM method for the 10th step and (d) is the results using the proposed method at the 10th step.
The regularization parameters for TV and first-order TV are selected manually as optimal as possible. 
The regularization parameters for the proposed method are set to be $\lambda=1 \times 10^{-10}, \rho=1 \times 10^{-9}, \delta = 0.06$ and $M=10$. For the comparison, in Fig. \ref{1D profiles} (a)-(c) we respectively plot the profiles along the solid, the dash and the dot dash lines respectively shown in Fig. \ref{numerical_results} (d). Fig. \ref{quantitative indicators} (a)-(b) show the behavior of $RE(n)$ and $PSNR(n)$ as $n$ increases for the methods excluding FER. In Table II we compare the computational time for the reconstructions using different regularizers.

From the comparisons, we obtained similar conclusions as the 2D case does.

\begin{table}
\centering
\caption{Computational times of the the 7th model in Fig. \ref{lung_results} for 2D simulation and the 3D simulations.}
\begin{tabular}{ccccc}
\hline\hline
Dimension                & FER      & TV       & First-order TV & NWFOTV   \\ \hline
2D & 0.029 & 2.311 & 0.644       & 0.629 \\
3D & 0.028 & 1.733 & 0.398 & 0.428 \\
\hline\hline
\end{tabular}
\end{table}

\subsection{EIT lung imaging experiment}

\begin{figure}[h]
\centering
\subfloat[]{\label{EIT measurement}
\includegraphics[width=0.3\linewidth]{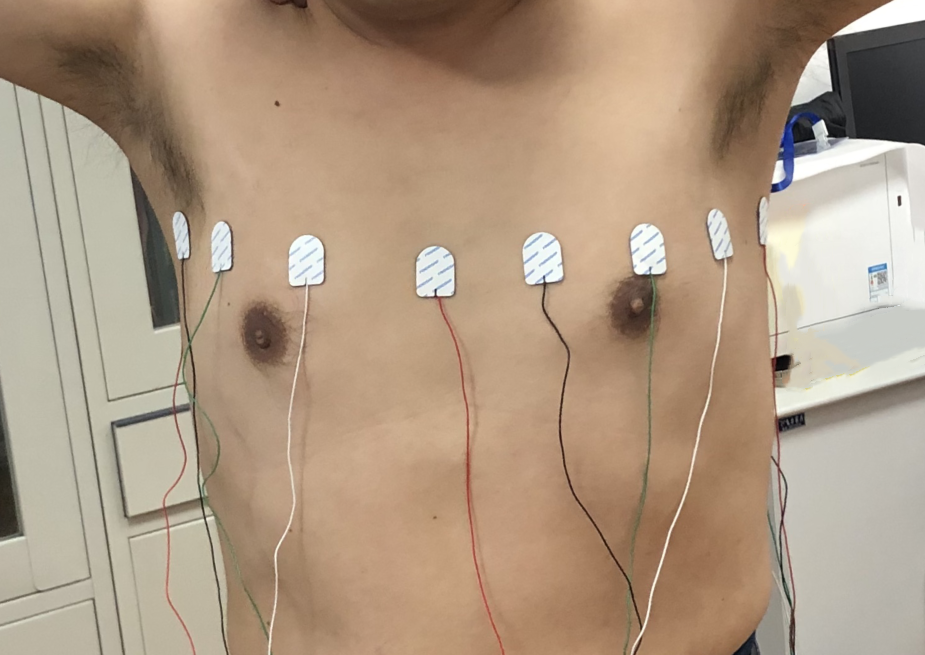}}
\subfloat[]{\label{3D_geometry}
\includegraphics[width=0.3\linewidth]{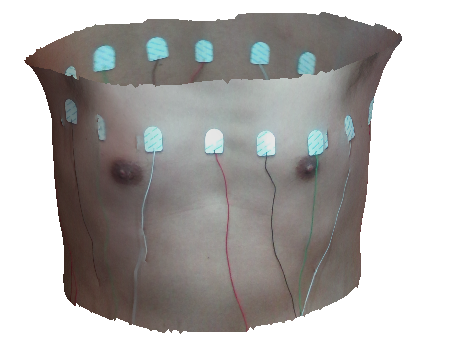}}
\subfloat[]{\label{imaging slice}
\includegraphics[width=0.3\linewidth]{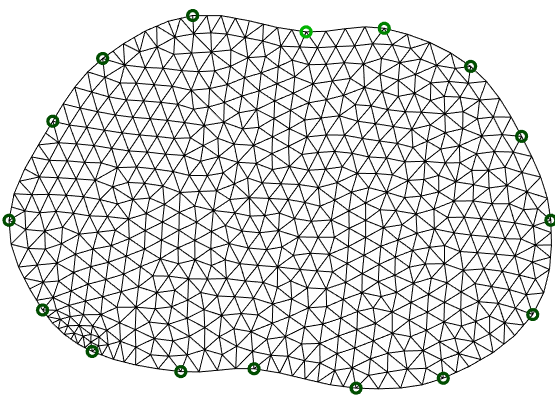}}
\centering
\caption{(a) Experiment setup, (b) the geometry obtained from 3D scanner, (c) a 2D imaging slice with point electrodes approximations.}
\label{Human experiment mesurement systerm}
\end{figure}

The human lung EIT experiment is approved by the ethics committee of the science and technology division, Shandong Normal University.
In this experiment, we attach 16 electrodes \cite{Company} around the exterior of the object's torso as shown in Fig. \ref{Human experiment mesurement systerm} (a). We choose the disposable Ag/AgCl electrocardiogram (ECG) electrodes with the size of $20\times 27$ mm. We use the 3D scanner \cite{Sense3D} to scan the body with electrodes and obtain the geometry and the positions of the electrodes as accurate as possible.
The geometry obtained is shown in Fig. \ref{Human experiment mesurement systerm} (b). To obtain the position of the electrodes, we use point electrode  approximations which lie on the center of the true electrodes \cite{Hanke2011}.
Then we discretize the imaging slice with point electrodes into 614 nodes and 1132 triangles. The discretization result is shown in Fig. \ref{Human experiment mesurement systerm} (c). During the reconstruction, we assign the reference conductivity $\sigma_0 = 0.0107$ S/m \cite{Adler2006} as the conductivity of the lungs at the state of the end expiration and use this $\sigma_0$ to obtain the sensitivity matrix $\mathbf{S}$ in \eqref{eq:sensitivity}.

\begin{figure*}[]
  \centering
  \includegraphics[width=0.9\linewidth]{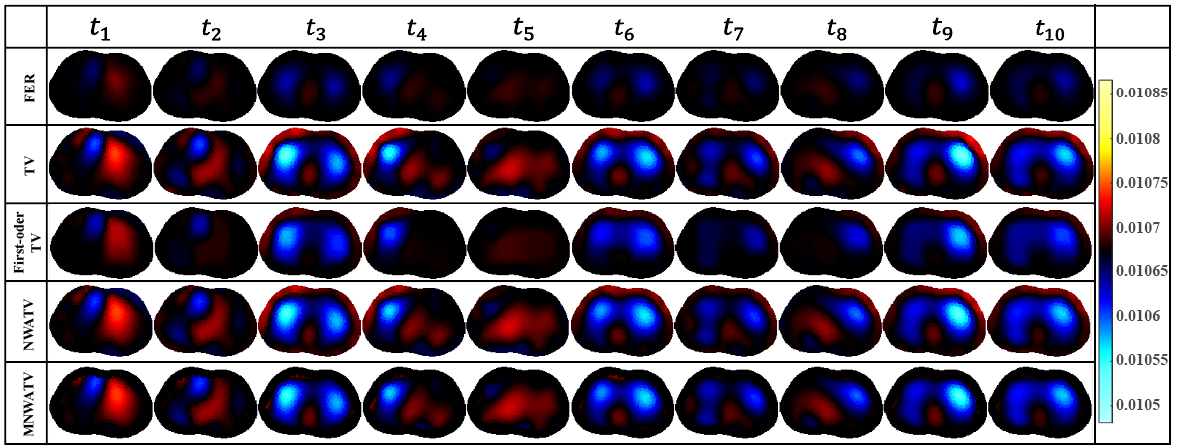}\\
  \caption{10 frames of human lung images using different regularizers. The first to the third row illustrate the reconstruction results using the FER, the TV and the first-order TV regularizers while the fourth to the fifth row illustrate the results using the NWATV and MNWATV regularizers respectively for $t=t_1,\cdots,t_{10}$.}\label{Human_experiment_Results}
\end{figure*}

The human EIT data was collected with the Sciospec 16-channel EIT system \cite{EIT}.
The frequency of the injection current is 10 kHz and the speed of the data acquisition is 30 frames/s.
Using the datum obtained from the EIT scanner, we reconstruct 10 frames of the human lung images.
The time step for each frame is around 0.66 s ($t_n - t_{n-1} \approx 0.66$ s for $n=1,2,\cdots,10$). 
In Fig. \ref{Human_experiment_Results} we compare the reconstructions using TV, FER, first-order TV, the NWATV and the modified NWATV regularizer (MNWATV). In this experiment the regularization parameters for the MNWATV regularizer are set to be $\lambda=1\times 10^{-14}$, $\rho=1.25 \times 10^{-5}$, $\delta=7\times 10^{-8}$, $M= 4$ and $\lambda_b=1 \times 10^{-7}$. Again for the regularizers of TV and the first order TV, the regularization parameters are set to be as optimal as we can.

Due to the noise in the measured datum, in order to block the error propagation in the iteration process, in each step we modify the reconstruction result by formula \eqref{modify_result}.
The results of the reconstructions using each regularizer are shown in Fig. \ref{Human_experiment_Results}. Especially the results of the modified reconstructions are shown in the last row of Fig. \ref{Human_experiment_Results}. Table \ref{human_RE_PSNR} illustrates the behavior of $\widetilde{RE}(n)$ and $\widetilde{PSNR}(n)$ for the frames $t_3$, $t_6$ and $t_{10}$. In Table IV we provide the computational time when using different regularizers.

As we can see through the use of the proposed method, we obtain a similar reconstruction as observed in TV images and a more accurate image than FER.
Moreover, the computational cost is acceptable in clinical applications since it produces up to 6 frames per seconds,
which is much higher than the typical respiratory rate.

\begin{table}[]
\centering
\caption{The behaviors of $\widetilde{RE}(n)$ and $\widetilde{PSNR}(n)$ for the frames $t_3$, $t_6$ and $t_{10}$}
\begin{tabular*}{8 cm}{@{\extracolsep{\fill}}cccccc}
\hline\hline
   & FER, $\lambda=\infty$    & First-order TV & MNWATV \\ \hline
   & \multicolumn{3}{c}{$\widetilde{RE}$}            \\ \cline{2-4}
$t_3$  & 0.0107 & 0.0039         & 0.0021  \\
$t_{6}$   & 0.0089 & 0.0037         & 0.0022  \\
$t_{10}$  & 0.0081 & 0.0034         & 0.0017  \\ \hline
   & \multicolumn{3}{c}{$\widetilde{PSNR}$}          \\ \cline{2-4}
$t_3$   & 19.13  & 27.83          & 33.70   \\
$t_{6}$   & 20.63  & 28.19          & 35.91   \\
$t_{10}$  & 21.36  & 28.71          & 35.00   \\ \hline\hline
\end{tabular*}
\label{human_RE_PSNR}
\end{table}

\begin{table}[]
\centering
\caption{Comparison of the five methods in terms of computational time
for human experiments.}
\begin{tabular*}{8 cm}{@{\extracolsep{\fill}}cccccc}
\hline\hline
    & \multicolumn{5}{c}{Computational time (s)}                                                              \\ \cline{2-6}
 &      TV    & FER   & First-order TV   &  NWATV    & MNWATV \\ \hline
$t_1$    & 2.055 & 0.020 & 0.116 & 0.109  & 0.153   \\
$t_2$    & 2.131 & 0.017 & 0.104 & 0.105  & 0.140   \\
$t_3$    & 2.056 & 0.019 & 0.119 & 0.108  & 0.169   \\
$t_4$    & 2.071 & 0.018 & 0.109 & 0.113  & 0.163   \\
$t_5$    & 2.081 & 0.021 & 0.116 & 0.112  & 0.149   \\
$t_6$ & 2.088 & 0.022 & 0.107 & 0.129  & 0.158   \\
$t_7$ & 2.140 & 0.020 & 0.105 & 0.122  & 0.199   \\
$t_8$ & 2.136 & 0.019 & 0.110 & 0.109  & 0.176   \\
$t_9$ & 2.157 & 0.02  & 0.107 & 0.108  & 0.135   \\
$t_{10}$ & 2.161 & 0.019 & 0.112 & 0.1149 & 0.135   \\ \hline\hline
\end{tabular*}
\label{time_consumption}
\end{table}

\section{Discussion: Regularization Parameter Selection and Future Work}\label{section:discussion}

TV regularization has the advantage of preserving the inter-medium discontinuities \cite{Javaheriana2016} especially for a piecewise constant images and definitely has promising in EIT clinical applications. \cite{Dobson1994} introduced the TV regularization in EIT, \cite{Borsic2010} introduce a primal dual interior-point-method for the minimizing process. However, this method is time consuming and consequently it is difficult to apply in real time applications.
The proposed nonlinear weighted anisotropic TV regularization method reduces the computational time by using a weighted anisotropic TV regularizer. Moreover, this method is promising in producing a better result than using only traditional TV (isotropic TV) or anisotropic TV regularization. This is because the weight $1/|\na\sigma|^2$ could avoid possible pseudo-edges in isotropic TV or distortion along the coordinates axes in anisotropic TV regularizations \cite{Gonzalez2017}.

\begin{figure}[h]
\centering
\subfloat[]{
\includegraphics[width=0.48\linewidth]{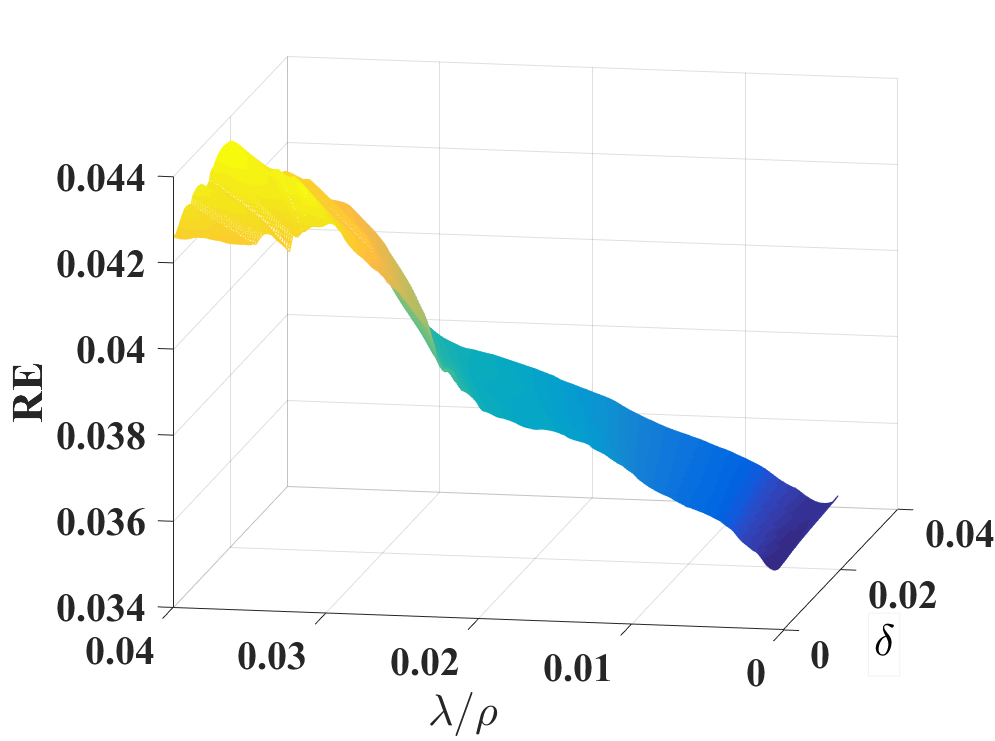}}
\subfloat[]{
\includegraphics[width=0.48\linewidth]{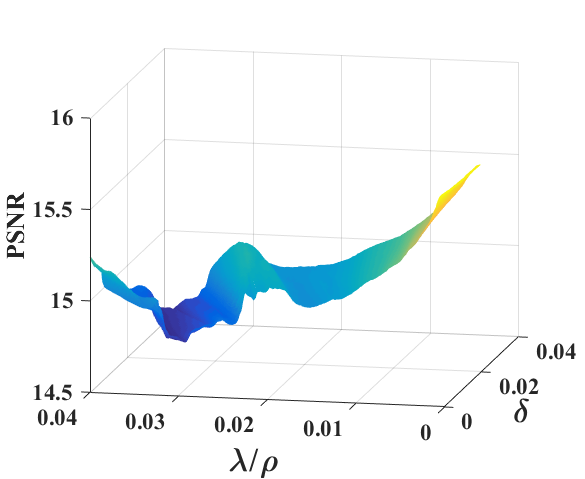}}
\centering
\caption{Changes of $RE$ and $PSNR$ with respect to $\frac{\lambda}{\rho}$ and $\delta$ in the 10th image of two-dimensional numerical simulation.}
\label{Change rule_of_quantitative indicators_in Simulation}
\end{figure}

For regularization based reconstructions, the choices of the regularization parameter(s) are always critical. In this paper, the choices of $\rho$ and $\delta$ depend on the mesh size used in the process of inversion. $\rho$ has been chosen in such a way that $\f{1}{\rho}\S^T \S$ and $\D^T \D$ in the similar order so that the formula \eqref{a51} makes sense. $\delta$ is chosen in the similar way as $\rho$ such that $\delta$ is in the similar order as $|\D\bmsigma|$. This explains why the parameter $\rho$ and $\delta$ is chosen so small. From the formula \eqref{a62}, only the ratio $\lambda/\rho$ is meaningful. Hence, the parameter $\lambda$ has to be chosen according to the selection of $\rho$.
Fig \ref{Change rule_of_quantitative indicators_in Simulation} (a) and (b) illustrate the change of $RE$ and $PSNR$ with respect to $\delta$ and $\lambda/\rho$, respectively. As we can see from this figure, for a fixed $\lambda/\rho$, $RE$ and $PSNR$ are almost invariant with respect to $\delta$. Moreover, there is optimal choice of $\lambda/\rho$ to minimize the $RE$ and maximize $PSNR$.

Future studies should cover a strict mathematical theory on the convergence analysis of the iteration process.
For clinical applications, a potential consideration is to add the low rank constraints of the sequential images \cite{Wang2019a} in the reconstruction.

\section{Conclusion}

In this paper, we propose a nonlinear weighted anisotropic TV regularization method in EIT to reduce the computational time and improve the capability of edge preservation in comparison with TV regularization based imaging. To validate the advantages of the proposed regularization method corroborating with the well established FER, TV, and first-order TV regularizers, we carried out 2D, 3D simulations and human EIT lung imaging.
In the testing campaign, it was shown that the proposed method reduces the computational time significantly while providing the reconstruction images similar to or even better than the traditional TV method.

\appendix
In this section we provide the proofs of the formula \eqref{a62}.

\subsection{Proof of the formula \eqref{a62}}

Direct calculation yields that
\begin{equation*}
\begin{split}
 \z_{n+1}& = \argmin_\z\lambda\|\mathbf{p}_n\cdot \mathbf{z}\|_{{l_1}} \\
 & \qquad\qquad +\frac{\rho}{2}\|\mathbf{D}\delta\bm{\sigma}_{n+1} - \mathbf{z}\|_{l_2}^{2}+\mathbf{y}_n^{T}(\mathbf{D}\delta\bm{\sigma}_{n+1} - \mathbf{z}) \\
& =\argmin_\z \|\mathbf{D}\delta\bm{\sigma}_{n+1} - \mathbf{z}\|_{l_2}^{2} \\
& \qquad\qquad +\f{2}{\rho}\mathbf{y}_n^{T}(\mathbf{D}\delta\bm{\sigma}_{n+1} - \mathbf{z}) + 2\f{\lambda}{\rho}\|\mathbf{p}_n\cdot \mathbf{z}\|_{{l_1}}\\
& = \argmin_\z\left\|\z-(\D\delta\bmsigma_{n+1}+\f{\y_n}{\rho})\right\|_{{l_2}}^2\\
&\qquad\qquad +\f{2\lambda}{\rho} \|\bfp_n\cdot \z\|_{l_1} -\left\|\f{\y_n}{\rho}\right\|_{{l_2}}^2 \\
& =\argmin_\z\left\|\z-(\D\delta\bmsigma_{n+1}+\f{\y_n}{\rho})\right\|_{{l_2}}^2+\f{2\lambda}{\rho} \|\bfp_n\cdot \z\|_{l_1},
\end{split}
\end{equation*}
where the last equality comes from the fact that $\|\y_n/\rho\|_{{l_2}}^2$ is independent of $\z$.
From the definition of $l_p$ ($p=1,2$) norm we obtain that
\begin{equation*}
  \z_{n+1} = \argmin\sum_{k=1}^N F_{\rho,\lambda}(\z[k]),
\end{equation*}
where $F_{\rho,\lambda}[\z[k]]$ is defined as follows
\begin{equation*}
  \begin{split}
    F_{\rho,\lambda}[\z[k]] & = \left\|\z[k]-\left(\D\delta\bmsigma_{n+1}[k]+\f{\y_n[k]}{\rho}\right)\right\|_{{l_2}}^2 \\
    &\qquad\qquad + \f{2\lambda}{\rho} |\bfp_n[k]\z[k]|.
  \end{split}
\end{equation*}
Hence
\begin{equation}\label{L1_minimizer1}
\begin{split}
  \z_{n+1}[k] = & \argmin F_{\rho,\lambda}(\z[k]).
\end{split}
\end{equation}

Next we will calculate the minimizer of $F_{\rho,\lambda}(\z[k])$, $\z_{n+1}[k]$ explicitly. Indeed, from the first-order optimal condition we obtain
\begin{equation*}
\begin{split}
  0& = \left.\f{1}{2}\f{\p F_{\rho,\lambda}(\z_{n+1}[k])}{\p z}\right|_{z = \z_{n+1}[k]} \\
  &= \z_{n+1}[k]-\left(\D\delta\bmsigma[k]+\f{\y[k]}{\rho}\right) + \f{\lambda\bfp[k]}{\rho}\sgn{\z_{n+1}[k]}.
\end{split}
\end{equation*}
Hence $\z_{n+1}[k]$ satisfies the following identity
\begin{equation}\label{relation:z_and_Q}
  \z_{n+1}[k] = \D\delta\bmsigma_{n+1}[k]+\f{\y_n[k]}{\rho} - \f{\lambda\bfp_n[k]}{\rho}\sgn{\z_{n+1}[k]}.
\end{equation}

From \eqref{relation:z_and_Q} and the fact that $\lambda \bfp[k]/\rho\geq 0$ for $k=1,2,\cdots,N$, it is obvious that
\begin{equation}\label{sign_relation_z_and_Q}
\sgn{\z_{n+1}[k]} = \sgn{\D\delta\bmsigma_{n+1}[k]+\y_n[k]/\rho}.
\end{equation}

Noting that from the basics of calculus, the candidates of $\z_{n+1}[k]$ is either the stationary point or the non-differentiable point of $F_{\rho,\lambda}$. Hence

\begin{itemize}
\item For the case $|\D\delta\bmsigma_{n+1}[k]+\y_n[k]/\rho|>\f{\lambda \bfp_n[k]}{\rho}$.
\end{itemize}
From \eqref{def:h_g}, \eqref{relation:z_and_Q} and \eqref{sign_relation_z_and_Q} we obtain
\begin{equation*}
  \z_{n+1}[k] = (I - \lambda \bfp_n[k]/\rho\cdot\mbox{sgn})[\D\delta\bmsigma_{n+1}[k]+\y_n[k]/\rho]
\end{equation*}
or $\z_{n+1}[k] = 0$, where $I$ is the identity map.

Since
\begin{equation*}
\begin{split}
&~~~~F_{\rho,\lambda}(\D\delta\bmsigma_{n+1}[k]+\y_n[k]/\rho]\mp \lambda \bfp_n[k]/\rho) \\
&= \f{\lambda^2}{\rho^2}(\bfp_n[k])^2 + \f{2\lambda}{\rho}\bfp_n[k]\cdot \\
&\qquad\qquad\qquad\qquad \left[\pm\left(\D\delta\bmsigma_{n+1}[k]+\f{\y_n[k]}{\rho}\right) \mp\f{\lambda}{\rho}\bfp_n[k]\right] \\
& =
\left\{
  \begin{split}
    & \f{2\lambda}{\rho}\bfp_n[k]\left(\D\delta\bmsigma_{n+1}[k]+\f{\y_n[k]}{\rho}\right)-\f{\lambda^2}{\rho^2}\bfp_n^2[k],\\
    &\qquad\qquad\qquad \mbox{if } \D\delta\bmsigma_{n+1}[k]+\f{\y_n[k]}{\rho}>\f{\lambda\bfp_n[k]}{\rho}\\
    & -\f{2\lambda}{\rho}\bfp_n[k]\left(\D\delta\bmsigma_{n+1}[k]+\f{\y_n[k]}{\rho}\right)+3\f{\lambda^2}{\rho^2}\bfp_n^2[k],\\
    &\qquad\qquad\qquad \mbox{if } \D\delta\bmsigma_{n+1}[k]+\f{\y_n[k]}{\rho}<-\f{\lambda\bfp_n[k]}{\rho}\\
  \end{split}
\right. \\
&<\left\|\D\delta\bmsigma_{n+1}[k]+\f{\y_n[k]}{\rho}\right\|_{l_2}^2 = L(0).
\end{split}
\end{equation*}
Hence in this case, \eqref{sign_relation_z_and_Q} holds.

\begin{itemize}
\item For the case $|\D\delta\bmsigma_{n+1}[k]+\y_n[k]/\rho|\leq \f{\lambda \bfp_n[k]}{\rho}$.
\end{itemize}
In this case, if $\z_{n+1}[k]\neq 0$ we have the following relation
\begin{equation*}
\begin{split}
& \sgn{(I - \lambda \bfp_n[k]/\rho\mbox{sgn})[\D\delta\bmsigma_{n+1}[k]+\y_n[k]/\rho]}\\
&\qquad\qquad \qquad\qquad\qquad\qquad\qquad\quad+\sgn{\z_{n+1}[k]}\equiv 0.
\end{split}
\end{equation*}
This contradicts with the relation \eqref{sign_relation_z_and_Q}. Hence the only minimizer of \eqref{L1_minimizer1} is $\z_{n+1}[k] = 0$.

Finally, we obtain the formula \eqref{a62}.

\end{document}